\documentclass[a4paper]{jpconf}
\usepackage{graphicx}

\usepackage{amsmath,amssymb,cite,comment}

\usepackage{amsfonts,stmaryrd}

\def\beq{\begin{equation}}
\def\eqn#1{\beq\label{#1}}
\def\ee{\end{equation}}
\def\eeq{\end{equation}}

\def\bb {\begin {eqnarray}}
\def\eqnn#1{\bb\label{#1}}
\def\eea {\end {eqnarray}}

\newcommand{\eqna}[1]{\begin{subequations} \label{#1}
\begin{eqnarray}}
\def\eena{\end{eqnarray}
\end{subequations}}

\def\th{{\tilde h}}

\def\hh{{\hat h}}

\def\nn{\nonumber}

\def\nd{\end{document}}

\input epsf.tex
\newcount\figno
\def\fig#1#2#3{
\par\begingroup\parindent=0pt\leftskip=1cm\rightskip=1cm\parindent=0pt
\baselineskip=11pt \global\advance\figno by 1 
\epsfxsize=#3 \centerline{\epsfbox{#2}} \vskip 12pt
#1\par
\endgroup\par}
\def\figlabel#1{\xdef#1{\the\figno}}
\def\encadremath#1{\vbox{\hrule\hbox{\vrule\kern8pt\vbox{\kern8pt
\hbox{$\displaystyle #1$}\kern8pt} \kern8pt\vrule}\hrule}}

  \def\tV{{\tilde V}}

\def\riga{-\kern-4pt - \kern-4pt -}
\font\fat=cmsy10 scaled\magstep5

\def\Bbullet{\raise-3pt\hbox{\fat\char"0F}}

\font\tfont=cmbx12 scaled\magstep1 

\def\Box{
\vbox{ \halign to5pt{\strut##& \hfil ## \hfil \cr &$\kern -0.5pt
\sqcap$ \cr \noalign{\kern -5pt \hrule} }}~}

\def\down{\raise1.5pt\hbox{$\phantom{a}_2$}\downarrow}

\def\downa{\raise1.5pt\hbox{$\phantom{a}_{2\atop m_2}$}\downarrow}

\def\llr{\longrightarrow}

\def\({\left(}
\def\){\right)}

\def\lra{\longrightarrow}

\def\bbc{\mathbb{C}}
\def\bac{\mathbb{C}}

\def\bbr{\mathbb{R}}
\def\bbn{\mathbb{N}}

\def\a{\alpha}
\def\b{\beta}
\def\d{\delta}

\def\vr{\vert}

\def\l{\lambda}

\def\D{{\Delta}}

\def\ca{{\cal A}}  \def\cc{{\cal C}}
\def\cd{{\cal D}}  \def\cf{{\cal F}}
\def\cg{{\cal G}} \def\ch{{\cal H}} 
 \def\ck{{\cal K}} 
\def\cm{{\cal M}} \def\cn{{\cal N}} 
\def\cp{{\cal P}}  
 \def\ct{{\cal T}}

\def\ido{intertwining differential operator}
\def\idos{intertwining differential operators}

\def\L{\Lambda}
\def\r{\rho}

\begin{document}

\begin{center}

{\tfont Langlands Duality and Invariant Differential Operators:\\[3pt]
the Case SL(2n+1)}

\vskip 1.5cm

{\bf V.K. Dobrev}

 \vskip 5mm

  Institute for Nuclear Research and Nuclear Energy,\\ Bulgarian
Academy of Sciences,\\ 72 Tsarigradsko Chaussee,  1784 Sofia,
Bulgaria

\end{center}

\vskip 1.5cm

 \centerline{{\bf Abstract}}

Recently we started building a bridge between two cases of Langlands duality. The latter is one of the most influential topics in mathematical research.  It has many different appearances and influential subtopics. Yet there is a topic that until now seems unrelated to the Langlands program.  That is the topic of invariant differential operators. 
 It is strange since both items are deeply rooted in Harish-Chandra’s represen
tation theory of semisimple Lie groups. Earlier, we studied the case of the
group SL(2n,R). In the present paper we deal with the group SL(2n+1,R).
The two groups just mentioned are similar, but their representation theories
are rather different.


\vskip 1.5cm

\section{Introduction}

In the last 50 years Langlands duality is one of the most influential topics in mathematical research \cite{Lan1,Lan2}.  It has many different appearances and influential subtopics, cf. an incomplete list 	in \cite{VKD-L}, and for some later papers we refer to
\cite{BSV,ChNa,ChVe,Chua,DHKM,DiTe,Esp,GaTe,GaWi,HLM,Hove,Ike,Ima,JLN,KiNo,LZZ,MaOp,Mat,Nak,Suz,Col,Sch,
Nima,WangWei,GLR,DEG,HFu,PSc,DHKZ,CJSL,CHS,JaMo,DHXZ,MOi,DYa,KBa,XuYi,THove,NaYu,MSY,VWa,DaFi,Toda,JYa,ACEST,CSW,GSh,SPa,THC,BHo,DPe,KaPo,LSW,IkRa,BaDa,Lia,Gle}.
Note that some papers are written by authors who have created influential topics themselves. The last fact
stresses the omnipresence of the Langlands program.

Our attempt is based on our approach to the construction of invariant differential operators - for an exposition we refer to \cite{VKD1} which is based also on many papers. 
Our approach is deeply related to the Langlands general classification of  representations of real semisimple groups $G$ \cite{Lan2} taking into account the refinement by Knapp-Zuckermann \cite{KnZu}. Thus, the main ingredient in \cite{Lan2} are the parabolic subgroups $P=MAN$, such that $M$ is semisimple subgroup of our group $G$ under study, $A$ is abelian subgroup, $N$ is nilpotent subgroup preserved by the action $A$. Altogether, there is a local (Bruhat) decomposition of $G$ using   a subgroup $G'= P\tilde{N}$, where $\tilde{N}$ is a  nilpotent subgroup of $G$ isomorphic to $N$ also preserved by the action $A$, so that $G'$ is dense in $G$.
According to the construction of Langlands-Knapp-Zuckermann every admissible irreducible representation  of $G$ may be obtained as a subrepresentation of representations of $G$ induced by a representations of some $P$ (some class is enough - see details below).

Further, the present paper is organized as follows. In the next section we give a synopsis of our approach.
Then we apply this to the group ~$SL(2n+1,\bbr)$, using the Langlands duality of the subgroup $M$ used in the example.
The cases $n=1,2,3$ are exposed in separate sections.


\section{Synopsis: Canonical construction of invariant differential operators}

Let $G$ be a semi-simple, non-compact Lie group, and $K$ a
maximal compact subgroup of $G$. Then, we have an {\it Iwasawa
decomposition} $G=KA_0N_0$, where $A_0$ is an Abelian simply
connected vector subgroup of $G$ and $N_0$ is a nilpotent simply
connected subgroup of $G$ preserved by the action of $A_0$.
Furthermore, let $M_0$ be the centralizer of $A_0$ in $K$. Then, the
subgroup $P_0 = M_0 A_0 N_0$ is a {\it minimal parabolic subgroup} of
$G$.

Further, for simplicity, we  restrict to {\it maximal parabolic
subalgebras} $\cp = \cm \oplus \ca \oplus \cn$ with $\dim\, \ca=1$.

Let $\nu$ be a (non-unitary) character of $A$, $\nu\in\ca^*$,
parameterized by a real number {\it $d$}, called  the {\it conformal weight}.
Next, $\mu$ fix a discrete series representation
 $D^\mu$ of $M$ on the Hilbert space $V_\mu\,$, or   the
finite-dimensional (non-unitary) representation of $M$ with the same
Casimirs.

We consider the induced
representation $\chi =$ Ind$^G_{P}(\mu\otimes\nu \otimes 1)$ an
 {\it \it elementary representation} of $G$ \cite{DMPPT},  (or   {\it generalized principal series representations} (or {\it
limits thereof}) in \cite{Knapp}.)  Their spaces of functions are  \begin{equation}\label{func}
\cc_\chi = \{ \cf \in C^\infty(G,V_\mu) \vr \cf (gman) =
e^{-\nu(H)} \cdot D^\mu(m^{-1})\, \cf (g) \} \end{equation} where $a=
\exp(H)\in A'$,~$H\in\ca'\,$,~$m\in M'$, $n\in N'$. The
representation action is the {\it left regular action}:  \begin{equation}\label{lrega}
(\ct^\chi(g)\cf) (g') = \cf (g^{-1}g'), \quad g,g'\in G.\end{equation}

Further we need   the {\it \it
highest/lowest-weight representations} of $\cg^\bac$. These can be
realized as (factor-modules of) Verma modules~$V^\L$ over
 $\cg^\bac$, where $\L\in (\ch^\bac)^*$, $\ch^\bac$ is a Cartan
subalgebra of $\cg^\bac$ and the weight $\L = \L(\chi)$ is determined
uniquely from $\chi$ \cite{VKD1}.

It is important to note that   Verma modules are
always reducible. Thus,  we use~{\it \it
generalized Verma modules}~$\tV^\L$~such that the role of the
highest/lowest-weight vector $v_0$ is taken by the
(finite-dimensional) space~$V_\mu\,v_0\,$. For the generalized
Verma modules (GVMs) the reducibility is controlled only by the
value of the conformal weight $d$. Relatedly, for the \idos{}, only
the reducibility with regard to non-compact roots is~essential.

 {Another}   main ingredient of our approach is as follows. We group
the (reducible) ERs with the same Casimirs in sets called {\it
multiplets} \cite{VKD1}. The multiplet corresponding to fixed values of the
Casimirs may be depicted as a connected graph, the {\it vertices} of
which correspond to the reducible ERs and the {\it lines (arrows)}
between the vertices correspond to intertwining operators.  The
explicit parameterization of the multiplets and of their ERs is
important in understanding of the situation. The notion of multiplets was introduced in \cite{Dobmul}
and applied to representations of~$SO_o(p,q)$~and~$SU(2,2)$, resp.,
induced from their minimal parabolic subalgebras. Then it was applied to
the conformal superalgebra \cite{DoPemul}, to quantum groups, 
infinite-dimensional (super)algebras 
- see volumes following  \cite{VKD1}.

In fact, the multiplets contain explicitly all the data necessary to
construct the \idos{}. Actually, the data for each \ido{} consist
of the pair $(\b,m)$, where $\b$ is a (non-compact) positive root
of $\cg^\bac$, $m\in\bbn$, such that the {\it BGG  Verma module
reducibility condition} (for highest-weight modules) is fulfilled:\vspace{-3pt}
\begin{equation}\label{bggr} (\L+\r, \b^\vee ) = m, \quad \b^\vee \equiv 2 \b
/(\b,\b) \ \end{equation}
where $\r$ is half the sum of the positive roots of
 $\cg^\bac$. When the above holds, then the Verma module with shifted
weight $V^{\L-m\b}$ (or $\tV^{\L-m\b}$ for GVM and $\b$
non-compact) is embedded in the Verma module $V^{\L}$ (or
 $\tV^{\L}$). This embedding is realized by a singular vector
 $v_s$ determined by a polynomial $\cp_{m,\b}(\cg^-)$ in the
universal enveloping algebra $(U(\cg_-))\ v_0\,$, and $\cg^-$ is the
subalgebra of $\cg^\bac$ generated by the negative root generators
\cite{Dix}.
 More explicitly \cite{VKD1}, $v^s_{m,\b} = \cp_{m,\b}\, v_0$ (or $v^s_{m,\b} = \cp_{m,\b}\, V_\mu\,v_0$ for GVMs).
   Then,
there exists \cite{VKD1}
 {\bf \ido{}} \begin{equation}\label{invop}  \cd_{m,\b} : \cc_{\chi(\L)}
 \llr \cc_{\chi(\L-m\b)} \end{equation} given explicitly by: \begin{equation}\label{singvv}
 \cd_{m,\b} = \cp_{m,\b}(\widehat{\cg^-})  \end{equation} where
 $\widehat{\cg^-}$ denotes the {\it right action} on the functions
 $\cf$.


In our exposition below, we shall use the so-called Dynkin labels:  \begin{equation}\label{dynk} m_i
 \equiv (\L+\r,\a^\vee_i) , \quad i=1,\ldots,n, \end{equation} where $\L =
\L(\chi)$, $\r$ is half the sum of the positive roots of
 $\cg^\bac$.

We shall use also   the so-called ~{\it Harish--Chandra parameters}:
\begin{equation}\label{dynhc} m_\b \equiv (\L+\r, \b^\vee )\ ,  \end{equation} where $\b$ is any
positive root of $\cg^\bac$. These parameters are redundant, since
they are expressed in terms of the Dynkin labels; however,   some
statements are best formulated in their terms. (Clearly, both the Dynkin labels and
Harish--Chandra parameters have their origin in the BGG reducibility condition \eqref{bggr}.)

 Next, we recall the action of the Weyl group on highest weights:
\eqn{weeyl}
w_\b(\L) ~\doteq~ \L - (\L+\r,\b^\vee) \b \ee
and thus,
\eqn{wezyl} w_\b(\L) ~=~ \L - m_\b \b \ee
and the shifted weight in \eqref{invop} results by the action of the Weyl group as in \eqref{wezyl}.

\section{The case of $SL(2n+1,\bbr)$}

In this paper we treat the case of $G=SL(2n+1,\bbr)$, $\cg=sl(2n+1,\bbr)$. We restrict to maximal parabolic subalgebra
\eqnn{parab} \cp ~&=&~ \cm \oplus \ca \oplus \cn \\
\cm &=& sl(n+1,\bbr) \oplus sl(n,\bbr) \,, ~\dim \ca = 1, ~\dim \cn = n^2 +n\eea

 In this section we start with $G=SL(m,\bbr)$, the group of invertible $m \times m $ matrices with real elements and
determinant 1. Then $\cg = sl(m,\bbr)$ and the Cartan involution is given explicitly by:
~$\th X = -\ ^tX$, where $^tX$ is the transpose of $X\in\cg$. Thus, $\ck \cong so(m,\bbr)=so(m)$, and is spanned by matrices
(r.l.s. stands for real linear span):
\eqn{maxc}  \ck = {\rm r.l.s.} \{ X_{ij} \equiv e_{ij} - e_{ji} \ , \quad 1\leq i < j \leq m \} \ ,\ee
where $e_{ij}$ are the standard matrices with only nonzero entry (=1) on the $i$-th row and $j$-th column,
$(e_{ij})_{k\ell} = \d_{ik}\d_{j\ell}\,$. (Note that $\cg$ does not have discrete series representations
if $m>2$.)

Further, the complementary space $\cp$ is given by:
\eqnn{slp} \cp &=& {\rm r.l.s.} \{ Y_{ij} \equiv e_{ij} + e_{ji} \ , ~~ 1\leq i < j \leq m \ ,\\
&&  H_j \equiv e_{jj} - e_{j+1,j+1} \ , \quad 1\leq j \leq m-1 \}\ .\eea
The split rank is ~$r = m-1$, and from \eqref{slp} it is obvious that in this setting one has:
\eqn{aaa}\ch = {\rm r.l.s.} \{ H_j   \ , \quad 1\leq j \leq n-1=r \}\ .\ee

The simple root vectors are given explicitly by:
\eqn{simple} X^+_j \doteq e_{j,j+1} \ , ~~~ X^-_j \doteq e_{j+1,j} \ , ~~~ 1\leq j \leq m-1 \ .\ee
Note that matters are arranged so that
\eqn{slsub} [X^+_j,X^-_j] = H_j\ , \quad [H_j,X^\pm_j] = \pm 2X^\pm_j\ ,\ee
and further we shall denote by ~$sl(2,\bbr)_j$~ the ~$sl(2,\bbr)$~
subalgebra of $\cg$ spanned  by $X^\pm_j\,, H_j\,$.

In our case of consideration ~$m=2n+1$ we have
\eqn{facm}  \cm ~=~ sl(n+1) \oplus sl(n) \ee
and we use representations of $\cm$ indexed as follows:
\eqn{repm} \hat{\cm} = (m_1,\ldots,m_{n} \ ;\ m_{n+2}, \ldots, m_{2n}) \ee
When all $m_j$ are natural numbers   $\hat{\cm}$ indexes the unitary finite-dimensional irreps of $\cm$.

Next we recall that the number of ERs in the corresponding multiplets is equal
 to the following ratio of numbers of elements of Weyl groups:
\eqn{multi} \vr W(\cg^\bac,\ch^\bac)\vr\, /\, \vr
W(\cm^\bac,\ch_m^\bac)\vr   \ee
where ~$\ch^\bac,\ \ch^\bac_m$~ are Cartan subalgebras of ~$\cg^\bac,\ \cm^\bac$, resp.

  The same number holds for
any algebra $\cg'$ parabolically related to $\cg$ w.r.t. $\cm$.

 Further, we denote by ~$\cc^\pm_i$~ the representation space with signature
~$\chi^\pm_i\,$.

We need to take into account the operators intertwining the pairs ~$\cc^\pm_i$~:
\eqn{knapps}
G^\pm_i ~:~ \cc^\mp_i \lra \cc^\pm_{i}   \ , \quad i ~=~
1,\ldots,1+\hh  \ . \ee
These are integral operators with kernels being described by the two-point functions that were used for
the first time in conformal field theories \cite{VKD1}.

{\bf Remark 1.} In the next sections we give explicitly the cases $sl(3),sl(4),sl(5)$. It is natural that we restrict to these cases for the reason of practicality. We recall that the numbers
of representations in a multiplet is given as follows:
\cite{VKD1}:
\eqn{redun} N_M ~=~  {\vert W(\cg,\ch) \vert \over \vert W(\cm,\ch_m) \vert} \ee
Thus, for $sl(3),sl(4),sl(5)$ the number of of ERs in the main multiplet is:
$$ N_M ~=~ 3, 10, 35 $$
Now if we go to the next case $sl(7)$ we shall have 126 representations in the multiplet.  Clearly, this is difficult to draw and also difficult to understand and interpret.
This is the reason of the restriction of the cases that we consider.

\section{The case of $SL(3,\bbr)$}

In the case of $sl(3)$ the  parabolic $\cm$ factor is:
\eqn{su22m} \cm_3 ~=~ sl(2) \ee
the representations being indexed by the number $m_1$.

When we consider   induction from $\cm_3 ~=~ sl(2) $ then we have  three-member multiplets
using   \
which in our case ($\cm=\cm_4$) gives:
\eqn{redu22} N_M ~=~  {6 \over 2} ~=~ 3. \ee

The multiplet is parametrized as follows:
\eqnn{tabsu22}
\chi^1 ~&=&~ \{\, m_1\,, m_2\,
\}  , \\
\chi^2 ~&=&~ \{\, m_{12}\,, -m_{2}\,
\}  , \ \L^2 = \L^1 - m_2\a_2 \nn\\
 \chi^3 ~&=&~ \{\, m_{2}\,, -m_{12}\,
\}  , \ \L^3  = \L^2 - m_1\a_{12} \nn\eea
 where  $m_{12}\equiv m_1+m_2$.

 It is well-known that when all $m_j$ are natural numbers then $\chi_1$ exhausts the finite-dimensional representations of $\cg$.

  Note that $\chi^1$ and $\chi^3$  are related by Knapp-Stein \cite{KnSt} integral intertwining operators ~$G^\pm$~ so that
the operators $G^+$ act from $\chi^1$ to $\chi^3$, while $G^-$ act from $\chi^3$ to $\chi^1$.

 Thus, the ~{\it Knapp-Stein duality}~ is a manifestation of the ~{\it Langlands duality}.

{\bf Remark: 2.}~ Note that unlike the cases of SL(2n) treated in \cite{VKD-L} there is no symmetry expressed through the conformal weight ~$d$.

Using the language of Weyl reflections we have:
\eqnn{slthr} && w_{\a_2} (\L^1_3) = \L^2_3 \\
&& w_{\a_{12}} (\L^1_3) = \L^2_3 \nn
\eea
while for the invariant differential operators 
correspondingly we  write:
\eqnn{slthr}  \cd_{m_2,\a_2} : \cc^1_3 \lra \cc^2_3 \\
\cd_{m_1,\a_{12}} : \cc^2_3 \lra \cc^3_3 \nn\eea

Diagramatically we have on Fig.1. the following $SL_3$ quiver:
\eqnn{qui3}Caption: SL_3\ quiver \nn\\ 
 {\bf Fig.1.} ~~~\L^1_3 \Rightarrow \L^2_3 \Rightarrow \L^3_3 \nn\eea

Multiplets containing the finite-dimensional subrepresentations of $\cg$ are called {\it main multiplets}.
 The other multiplets are called reduced multiplets.
 In the case of reduced quiver (multiplet) we have:
\eqnn{slthrc}
\chi^{'1}_3 &=& \{ m_1, 0 \   \} \ , ~~m_1 \in\bbn \\
 \chi^{'2}_3 &=& \{ 0, -m_{1}\ \}, ~~~\L^{'2}_3 = \L^{'1}_3 - m_1 \a_{12} \nn\eea

And instead of  \eqref{slthr} we have:
 \eqn{slthre} \cd_{m_1,\a_{12}} : \cc^{'1}_3 \lra \cc^{'2}_3 \ee

 It is important that the differential operator in \eqref{slthre} is obtained from the Knapp-Stein integral intertwining  operator after regularization
 of the integral kernel using the old mechanism described many years ago by Gelfand et al \cite{GeV}. After the regularization at these singular representation
 points the integral kernel produces a delta function which leads to the fact that the integral operator produces an \ido{}.

\section{The case of $SL(5,\bbr)$}

 The number of elements of the main multiplets   of a
  Lie algebra $\cg$ with $\cm$-factor fulfilling \eqref{repm}
is given by \eqref{multi}, i.e., 10 in the current case.


\eqnn{slt5}
\chi^1_5 &=& \{ m_1, m_2, m_3, m_4 \   \} \ , ~~m_j \in\bbn \\
\chi^2_5 &=& \{ m_1, m_{23}, -m_3, m_{34} \   \} , ~~~\L_2 = \L_1 - m_3 \a_3 \nn\\
\chi^3_5 &=& \{ m_{1}, m_{24}, -m_{34}, m_{3} \   \} , ~~~\L_3 = \L_2 - m_4 \a_{34} \nn\\
\chi^4_5 &=& \{ m_{12}, m_{3}, -m_{23}, m_{24} \   \} , ~~~\L_4 = \L_2 - m_2 \a_{23} \nn\\
\chi^5_5 &=& \{ m_{12}, m_{34}, -m_{24}, m_{23} \   \} , ~~~\L_5 = \L_3 - m_2 \a_{23} = \L_4 - m_4 \a_{34} \nn\\
\chi^6_5 &=& \{ m_{2}, m_{3}, -m_{13}, m_{14} \   \} , ~~~\L_6 = \L_4 - m_1 \a_{13} \nn\\
\chi^7_5 &=& \{ m_{13}, m_{4}, -m_{24}, m_{2} \   \} , ~~~\L_7 = \L_5 - m_3 \a_{24} \nn\\
\chi^8_5 &=& \{ m_{2}, m_{34}, -m_{14}, m_{13} \   \} , ~~~\L_8 = \L_5 - m_1 \a_{13} = \L_6 - m_4 \a_{34}\nn\\
\chi^9_5 &=& \{ m_{23}, m_{4}, -m_{14}, m_{12} \   \} , ~~~\L_9 = \L_7 - m_1 \a_{13} = \L_8 - m_3 \a_{24} \nn\\
\chi^{10}_5 &=& \{ m_{3}, m_{4}, -m_{14}, m_{1} \   \} , ~~~\L_{10} = \L_9 - m_2 \a_{14} \nn\eea

Diagramatically we have this quiver on Fig.2.

To present the Knapp-Stein duality, denoted here by $\sim$, we shall write:
\eqnn{slt5d} \chi^6_5 &\sim& \tilde{\chi}^5_5 \\
\chi^7_5 &\sim&  {\chi}^4_5 \nn\\
\chi^8_5 &\sim&  {\chi}^3_5 \nn\\
\chi^9_5 &\sim&  {\chi}^2_5 \nn\\
\chi^{10}_5 &\sim&  {\chi}^1_5 \nn\eea
{\bf Remark: 3.}~ And again we note that unlike the cases of SL(2n) treated in \cite{VKD-L} there is no symmetry expressed through a conformal weight.

{\bf Reduced multiplets (quivers)}

Subtype $m1$~:
\eqnn{slt51}
_1\chi^1_5 &=& \{ 0, m_2, m_3, m_4 \   \} \ ,  ~~m_j \in\bbn \\
_1\chi^2_5 &=& \{ 0, m_{23}, -m_3, m_{34} \   \} , ~~~\L_2  = \L_1 - m_3 \a_3 \nn\\
_1\chi^3_5 &=& \{ 0, m_{24}, -m_{34}, m_{3} \   \} , ~~~\L_3 = \L_2 - m_4 \a_{34} \nn\\
_1\chi^4_5 &=& \{ m_{2}, m_{3}, -m_{23}, m_{24} \   \} , ~~~\L_4 = \L_6 = \L_2 - m_2 \a_{23} \nn\\
_1\chi^5_5 &=& \{ m_{2}, m_{34}, -m_{24}, m_{23} \   \} , ~~~\L_5 = \L_8 = \L_3 - m_2 \a_{23} = \L_4 - m_4 \a_{34} \nn\\
_1\chi^7_5 &=& \{ m_{23}, m_{4}, -m_{24}, m_{2} \   \} ,  ~~~\L_7   = \L_9 = \L_5 - m_3 \a_{24} \nn\\
_1\chi^{10}_5 &=& \{ m_{3}, m_{4}, -m_{24}, 0 \   \} , ~~~\L_{10} = \L_7 - m_2 \a_{14} \nn\eea

Subtype $m2$~:
\eqnn{slt52}
_2\chi^1_5 &=& \{ m_1, 0, m_3, m_4 \   \} \ , ~~m_j \in\bbn  \\
_2\chi^2_5 &=& \{ m_1, m_{3}, -m_3, m_{34} \   \} , ~~~\L_2 = \L_4 = \L_1 - m_3 \a_3 \nn\\
_2\chi^3_5 &=& \{ m_{1}, m_{34}, -m_{34}, m_{3} \   \} , ~~~\L_3 = \L_5 = \L_2 - m_4 \a_{34} \nn\\
_2\chi^6_5 &=& \{ 0, m_{3}, -m_{1,3}, m_{1,34} \   \} , ~~~\L_6 = \L_4 - m_1 \a_{13} \nn\\
_2\chi^7_5 &=& \{ m_{1,3}, m_{4}, -m_{34}, 0 \   \} , ~~~\L_7 = \L_5 - m_3 \a_{24} \nn\\
_2\chi^8_5 &=& \{ 0, m_{34}, -m_{1,34}, m_{1,3} \   \} , ~~~\L_8 = \L_5 - m_1 \a_{13} = \L_6 - m_4 \a_{34}\nn\\
_2\chi^9_5 &=& \{ m_{3}, m_{4}, -m_{1,34}, m_{1} \   \} , ~~~\L_9 = \L_{10} = \L_7 - m_1 \a_{13} = \L_8 - m_3 \a_{24} \nn\\
\nn\eea

Subtype $m3$~:
\eqnn{slt53}
_3\chi^1_5 &=& \{ m_1, m_{2}, 0, m_{4} \   \} , ~~~\L_2 = \L_1, ~~m_j \in\bbn  \\
_3\chi^3_5 &=& \{ m_{1}, m_{2,4}, -m_{4}, 0 \   \} , ~~~\L_3 = \L_2 - m_4 \a_{34} \nn\\
_3\chi^4_5 &=& \{ m_{12}, 0, -m_{2}, m_{2,4} \   \} , ~~~\L_4 = \L_2 - m_2 \a_{23} \nn\\
_3\chi^5_5 &=& \{ m_{12}, m_{4}, -m_{2,4}, m_{2} \   \} , ~~~\L_5 = \L_7 =  \L_3 - m_2 \a_{23} = \L_4 - m_4 \a_{34} \nn\\
_3\chi^6_5 &=& \{ m_{2}, 0, -m_{12}, m_{12,4} \   \} , ~~~\L_6 = \L_4 - m_1 \a_{13} \nn\\
_3\chi^8_5 &=& \{ m_{2}, m_{4}, -m_{12,4}, m_{12} \   \} , ~~~\L_8 = \L_9 = \L_5 - m_1 \a_{13} = \L_6 - m_4 \a_{34}\nn\\
_3\chi^{10}_5 &=& \{ 0, m_{4}, -m_{142,}, m_{1} \   \} , ~~~\L_{10} = \L_9 - m_2 \a_{14} \nn\eea

Subtype $m4$~:
\eqnn{slt5}
_4\chi^1_5 &=& \{ m_1, m_2, m_3, 0 \   \} \ , ~~m_j \in\bbn \\
_4\chi^2_5 &=& \{ m_1, m_{23}, -m_3, m_{3} \   \} , ~~~\L_2 = \L_5 =\L_1 - m_3 \a_3 \nn\\
_4\chi^4_5 &=& \{ m_{12}, m_{3}, -m_{23}, m_{23} \   \} , ~~~\L_3 = \L_4 = \L_2 - m_2 \a_{23} \nn\\
_4\chi^6_5 &=& \{ m_{2}, m_{3}, -m_{13}, m_{13} \   \} , ~~~\L_6 = \L_8 = \L_4 - m_1 \a_{13} \nn\\
_4\chi^7_5 &=& \{ m_{13}, 0, -m_{23}, m_{2} \   \} , ~~~\L_7 = \L_5 - m_3 \a_{24} \nn\\
_4\chi^9_5 &=& \{ m_{23}, 0, -m_{13}, m_{12} \   \} , ~~~\L_9 = \L_7 - m_1 \a_{13} = \L_8 - m_3 \a_{24} \nn\\
_4\chi^{10}_5 &=& \{ m_{3},0, -m_{13}, m_{1} \   \} , ~~~\L_{10} = \L_9 - m_2 \a_{14} \nn\eea

For further reduced cases we show only those with physical applicability, i.e.,  unitary with respect $\cm$:
 \eqnn{slt513}
  _{14}\chi^4_5 &=& \{ m_{2}, m_{3}, -m_{23}, m_{23} \   \} , ~~~\L_4 = \L_6 = \L_2 - m_2 \a_{23} \nn\\
_{13}\chi^5_5 &=& \{ m_{2}, m_{4}, -m_{2,4}, m_{2} \   \} , ~~~\L_5 = \L_8 = \L_3 - m_2 \a_{23} = \L_4 - m_4 \a_{34} \nn\\
_{24}\chi^2_5 &=& \{ m_1, m_{3}, -m_3, m_{3} \   \} , ~~~\L_2 = \L_4 = \L_1 - m_3 \a_3 \nn
  \eea

\section{The case of $SL(7,\bbr)$}

 In this case the number of elements of the main multiplets   of a
  Lie algebra $\cg$ with $\cm$-factor fulfilling ~\eqref{repm} 
is given by \eqref{multi}, i.e., 35: 

\eqnn{slt7}
\chi^1_7 &=& \{ m_1, m_2, m_3, m_4, m_5, m_6 \   \} \ , ~~m_j \in\bbn \\
 \chi^2_7 &=& \{ m_1, m_2, m_{34}, -m_4, m_{45}, m_6 \   \} \ ,~~\L_2 = \L_1 - m_4 \a_4 \nn\\ 
 \chi^3_7 &=& \{ m_1, m_{23}, m_{4}, -m_{34}, m_{35}, m_6 \   \} \ ,~~\L_3 = \L_2 - m_3 \a_{34} \nn\\ 
 \chi^4_7 &=& \{ m_1, m_2, m_{34}, -m_{45}, m_{4}, m_{56} \   \} \ ,~~\L_4 = \L_2 - m_5 \a_{45} \nn\\ 
  \chi^5_7 &=& \{ m_{12}, m_{3}, m_{4}, -m_{24}, m_{25}, m_6 \   \} \ ,~~\L_5 = \L_3 - m_2 \a_{24} \nn\\ 
\chi^6_7 &=& \{ m_1, m_{23}, m_{45}, -m_{35}, m_{34}, m_{56} \   \} \ ,~~\L_6 = \L_3 - m_5 \a_{45} \nn\\ 
\chi^7_7 &=& \{ m_1, m_2, m_{36}, -m_{46}, m_{4}, m_5 \   \} \ ,~~\L_7 = \L_4 - m_6 \a_{46} \nn\\  
\chi^8_7 &=& \{ m_2, m_{3}, m_{4}, -m_{14}, m_{15}, m_6 \   \} \ ,~~\L_8 = \L_5 - m_1 \a_{14} \nn\\ 
\chi^9_7 &=& \{ m_{12}, m_{3}, m_{45}, -m_{25}, m_{24}, m_{56} \   \} \ ,~~\L_9 = \L_5 - m_5 \a_{45} \nn\\ 
\chi^{10}_7 &=& \{ m_1, m_{24}, m_{5}, -m_{35}, m_{3}, m_{46} \   \} \ ,~~\L_{10} = \L_6 - m_4 \a_{35}  \nn\\ 
\chi^{11}_7 &=& \{ m_1, m_{23}, m_{46}, -m_{36}, m_{34}, m_{5} \   \} \ ,~~\L_{11} = \L_6 - m_6 \a_{46}  \nn\\ 
\chi^{12}_7 &=& \{ m_2, m_{3}, m_{45}, -m_{15}, m_{14}, m_{56} \   \} \ ,~~\L_{12} = \L_8 - m_5 \a_{45} \nn\\ 
\chi^{13}_7 &=& \{ m_{12}, m_{3}, m_{46}, -m_{26}, m_{24}, m_{5} \   \} \ ,~~\L_{13} = \L_9 - m_6 \a_{46} \nn\\ 
\chi^{14}_7 &=& \{ m_{12}, m_{34}, m_{5}, -m_{25}, m_{23}, m_{46} \   \} \ ,~~\L_{14} = \L_9 - m_4 \a_{35} \nn\\ 
\chi^{15}_7 &=& \{ m_1, m_{24}, m_{56}, -m_{36}, m_{3}, m_{45} \   \} \ ,~~\L_{15} = \L_{10} - m_6 \a_{46}  \nn\\ 
\chi^{16}_7 &=& \{ m_2, m_{34}, m_{5}, -m_{15}, m_{13}, m_{46} \   \} \ ,~~\L_{16} = \L_{12} - m_4 \a_{35} \nn\\ 
\chi^{17}_7 &=& \{ m_{2}, m_{3}, m_{46}, -m_{16}, m_{14}, m_{5} \   \} \ ,~~\L_{17} = \L_{13} - m_1 \a_{14} \nn\\ 
\chi^{18}_7 &=& \{ m_{12}, m_{34}, m_{56}, -m_{26}, m_{23}, m_{45} \   \} \ ,~~\L_{18} = \L_{15} - m_2 \a_{24}  \nn\\ 
\chi^{19}_7 &=& \{ m_{13}, m_{4}, m_{5}, -m_{25}, m_{2}, m_{36} \   \} \ ,~~\L_{19} = \L_{14} - m_3 \a_{25} \nn\\ 
\chi^{20}_7 &=& \{ m_1, m_{25}, m_{6}, -m_{36}, m_{3}, m_{4} \   \} \ ,~~\L_{20} = \L_{15} - m_5 \a_{36}  \nn\\ 
\chi^{21}_7 &=& \{ m_{23}, m_{4}, m_{5}, -m_{15}, m_{12}, m_{36} \   \} \ ,~~\L_{21} = \L_{16} - m_3 \a_{25} \nn\\ 
\chi^{22}_7 &=& \{ m_{2}, m_{34}, m_{56}, -m_{16}, m_{13}, m_{45} \   \} \ ,~~\L_{22} = \L_{16} - m_6 \a_{46} \nn\\ 
\chi^{23}_7 &=& \{ m_{13}, m_{4}, m_{56}, -m_{26}, m_{2}, m_{35} \   \} \ ,~~\L_{23} = \L_{18} - m_3 \a_{25}  \nn\\ 
\chi^{24}_7 &=& \{ m_{12}, m_{35}, m_{6}, -m_{26}, m_{23}, m_{4} \   \} \ ,~~\L_{24} = \L_{18} - m_5 \a_{36}  \nn\\ 
\chi^{25}_7 &=& \{ m_{3}, m_{4}, m_{5}, -m_{15}, m_{1}, m_{26} \   \} \ ,~~\L_{25} = \L_{21} - m_2 \a_{15} \nn\\ 
\chi^{26}_7 &=& \{ m_{23}, m_{4}, m_{56}, -m_{16}, m_{12}, m_{35} \   \} \ ,~~\L_{26} = \L_{21} - m_6 \a_{46} \nn\\ 
\chi^{27}_7 &=& \{ m_{2}, m_{35}, m_{6}, -m_{16}, m_{13}, m_{4} \   \} \ ,~~\L_{27} = \L_{22} - m_5 \a_{36} \nn\\ 
\chi^{28}_7 &=& \{ m_{13}, m_{45}, m_{6}, -m_{26}, m_{2}, m_{34} \   \} \ ,~~\L_{28} = \L_{24} - m_3 \a_{25}  \nn\\ 
\chi^{29}_7 &=& \{ m_{3}, m_{4}, m_{56}, -m_{16}, m_{1}, m_{25} \   \} \ ,~~\L_{29} = \L_{25} - m_6 \a_{46} \nn\\ 
\chi^{30}_7 &=& \{ m_{23}, m_{45}, m_{6}, -m_{16}, m_{12}, m_{34} \   \} \ ,~~\L_{30} = \L_{26} - m_5 \a_{36} \nn\\ 
\chi^{31}_7 &=& \{ m_{14}, m_{5}, m_{6}, -m_{26}, m_{2}, m_{3} \   \} \ ,~~\L_{31} = \L_{28} - m_4 \a_{26}  \nn\\ 
\chi^{32}_7 &=& \{ m_{3}, m_{45}, m_{6}, -m_{16}, m_{1}, m_{24} \   \} \ ,~~\L_{32} = \L_{29} - m_5 \a_{36} \nn\\ 
\chi^{33}_7 &=& \{ m_{24}, m_{5}, m_{6}, -m_{16}, m_{12}, m_{3} \   \} \ ,~~\L_{33} = \L_{30} - m_4 \a_{26} \nn\\ 
\chi^{34}_7 &=& \{ m_{34}, m_{5}, m_{6}, -m_{16}, m_{1}, m_{23} \   \} \ ,~~\L_{34} = \L_{32} - m_4 \a_{26} \nn\\ 
\chi^{35}_7 &=& \{ m_{4}, m_{5}, m_{6}, -m_{16}, m_{1}, m_{2} \   \} \ ,~~\L_{35} = \L_{34} - m_3 \a_{16} \nn\eea 

Diagramatically we have this quiver on Fig.3.

To present the Knapp-Stein duality, denoted again by $\sim$, we shall write:
\eqnn{slt7d} \chi^{35}_7 &\sim&  {\chi}^1_7 \\
\chi^{34}_7 &\sim&  {\chi}^2_7 \nn\\
\chi^9_7 &\sim&  {\chi}^2_7 \nn\\
\chi^{16}_7 &\sim& \chi^{19}_7 \nn\\
\chi^{18}_7 &\sim& \chi^{20}_7 \nn\\
\chi^{8}_7 &\sim& \chi^{25}_7 \nn\\
\chi^{11}_7 &\sim& \chi^{27}_7 \nn\\
\chi^{9}_7 &\sim& \chi^{26}_7\nn\\
\chi^{10}_7 &\sim& \chi^{28}_7 \nn\\
 \chi^{12}_7 &\sim& \chi^{21}_7 \nn\\
\chi^{13}_7 &\sim& \chi^{22}_7 \nn\\
\chi^{14}_7 &\sim& \chi^{23}_7\nn\\
\chi^{15}_7 &\sim& \chi^{24}_7 \nn\eea

{\bf Remark: 4.}~ And again we note that unlike the cases of SL(2n) treated in \cite{VKD-L} there is no symmetry expressed through a conformal weight.

\subsection{Reduced multiplets}

\eqnn{slt71}
   _1\chi^5_7 &=& \{ m_{2}, m_{3}, m_{4}, -m_{24}, m_{25}, m_6 \   \} \ ,~ \\ 
 _1\chi^9_7 &=& \{ m_{2}, m_{3}, m_{45}, -m_{25}, m_{24}, m_{56} \   \} \ ,~~\L_9 = \L_{12} =\L_5 - m_5 \a_{45} \nn\\ 
  _1\chi^{13}_7 &=& \{ m_{2}, m_{3}, m_{46}, -m_{26}, m_{24}, m_{5} \   \} \ ,~~\L_{13} = \L_{17} = \L_9 - m_6 \a_{46} \nn\\ 
_1\chi^{14}_7 &=& \{ m_{2}, m_{34}, m_{5}, -m_{25}, m_{23}, m_{46} \   \} \ ,~~\L_{14} = \L_{16} =\L_9 - m_4 \a_{35} \nn\\ 
   _1\chi^{18}_7 &=& \{ m_{2}, m_{34}, m_{56}, -m_{26}, m_{23}, m_{45} \   \} \ ,~~\L_{18} = \L_{22} = \L_{15} - m_2 \a_{24}  \nn\\ 
_1\chi^{19}_7 &=& \{ m_{23}, m_{4}, m_{5}, -m_{25}, m_{2}, m_{36} \   \} \ ,~~\L_{19} = \L_{21} = \L_{14} - m_3 \a_{25} \nn\\ 
 _1\chi^{23}_7 &=& \{ m_{23}, m_{4}, m_{56}, -m_{26}, m_{2}, m_{35} \   \} \ ,~~\L_{23} = \L_{26} = \L_{18} - m_3 \a_{25}  \nn\\ 
_1\chi^{24}_7 &=& \{ m_{2}, m_{35}, m_{6}, -m_{26}, m_{23}, m_{4} \   \} \ ,~~\L_{24} = \L_{27} = \L_{18} - m_5 \a_{36}  \nn\\ 
   _1\chi^{28}_7 &=& \{ m_{23}, m_{45}, m_{6}, -m_{26}, m_{2}, m_{34} \   \} \ ,~~\L_{28} = \L_{30} = \L_{24} - m_3 \a_{25}  \nn\\ 
  _1\chi^{31}_7 &=& \{ m_{24}, m_{5}, m_{6}, -m_{26}, m_{2}, m_{3} \   \} \ ,~~\L_{31} = \L_{33} = \L_{28} - m_4 \a_{26}   
   \nn\eea 

\eqnn{slt72}
  _2\chi^3_7 &=& \{ m_1, m_{3}, m_{4}, -m_{34}, m_{35}, m_6 \   \} \ ,~~\L_3 = \L_3 = \L_2 - m_3 \a_{34} \\ 
    _2\chi^6_7 &=& \{ m_1, m_{3}, m_{45}, -m_{35}, m_{34}, m_{56} \   \} \ ,~~\L_6 = \L_9 = \L_3 - m_5 \a_{45} \nn\\ 
  _2\chi^{10}_7 &=& \{ m_1, m_{34}, m_{5}, -m_{35}, m_{3}, m_{46} \   \} \ ,~~\L_{10} = \L_{14} =\L_6 - m_4 \a_{35}  \nn\\ 
_2\chi^{11}_7 &=& \{ m_1, m_{3}, m_{46}, -m_{36}, m_{34}, m_{5} \   \} \ ,~~\L_{11} = \L_{13} = \L_6 - m_6 \a_{46}  \nn\\ 
   _2\chi^{15}_7 &=& \{ m_1, m_{34}, m_{56}, -m_{36}, m_{3}, m_{45} \   \} \ ,~~\L_{15} = \L_{18} = \L_{10} - m_6 \a_{46}  \nn\\ 
  _2\chi^{20}_7 &=& \{ m_1, m_{35}, m_{6}, -m_{36}, m_{3}, m_{4} \   \} \ ,~~\L_{20} = \L_{24} = \L_{15} - m_5 \a_{36}  \nn\\ 
_2\chi^{21}_7 &=& \{ m_{3}, m_{4}, m_{5}, -m_{1,35}, m_{1}, m_{36} \   \} \ ,~~\L_{21} = \L_{25} =\L_{16} - m_3 \a_{25} \nn\\ 
  _2\chi^{26}_7 &=& \{ m_{3}, m_{4}, m_{56}, -m_{36}, m_{1}, m_{35} \   \} \ ,~~\L_{26} = \L_{21} - m_6 \a_{46} \nn\\ 
  _2\chi^{29}_7 &=& \{ m_{3}, m_{4}, m_{56}, -m_{1,36}, m_{1}, m_{35} \   \} \ ,~~\L_{29} = \L_{25} - m_6 \a_{46} \nn\\ 
_2\chi^{30}_7 &=& \{ m_{3}, m_{45}, m_{6}, -m_{1,36}, m_{1}, m_{34} \   \} \ ,~~\L_{30} = \L_{32} =\L_{26} - m_5 \a_{36} \nn\\ 
 _2\chi^{33}_7 &=& \{ m_{34}, m_{5}, m_{6}, -m_{1,36}, m_{1}, m_{3} \   \} \ ,~~\L_{33} = \L_{34} = \L_{30} - m_4 \a_{26} \nn 
  \eea

  \eqnn{slt73}
 _3\chi^2_7 &=& \{ m_1, m_2, m_{4}, -m_4, m_{45}, m_6 \   \} \ ,~~\L_2 = \L_1 - m_4 \a_4 \\ 
  _3\chi^4_7 &=& \{ m_1, m_2, m_{4}, -m_{45}, m_{4}, m_{56} \   \} \ ,~~\L_4 = \L_2 - m_5 \a_{45} \nn\\ 
   _3\chi^6_7 &=& \{ m_1, m_{2}, m_{45}, -m_{45}, m_{4}, m_{56} \   \} \ ,~~\L_6 = \L_3 - m_5 \a_{45} \nn\\ 
_3\chi^7_7 &=& \{ m_1, m_2, m_{46}, -m_{46}, m_{4}, m_5 \   \} \ ,~~\L_7 = \L_{11} = \L_4 - m_6 \a_{46} \nn\\  
  _3\chi^{14}_7 &=& \{ m_{12}, m_{4}, m_{5}, -m_{2,45}, m_{2}, m_{46} \   \} \ ,~~\L_{14} = \L_{19} = \L_9 - m_4 \a_{35} \nn\\ 
 _3\chi^{16}_7 &=& \{ m_2, m_{4}, m_{5}, -m_{12,45}, m_{12}, m_{46} \   \} \ ,~~\L_{16} = \L_{21} = \L_{12} - m_4 \a_{35} \nn\\ 
  _3\chi^{18}_7 &=& \{ m_{12}, m_{4}, m_{56}, -m_{2,46}, m_{2}, m_{45} \   \} \ ,~~\L_{18} = \L_{23} = \L_{15} - m_2 \a_{24}  \nn\\ 
 _3\chi^{22}_7 &=& \{ m_{2}, m_{4}, m_{56}, -m_{12,46}, m_{12}, m_{45} \   \} \ ,~~\L_{22} = \L_{26} = \L_{16} - m_6 \a_{46} \nn\\ 
_3\chi^{24}_7 &=& \{ m_{12}, m_{45}, m_{6}, -m_{2,46}, m_{2}, m_{4} \   \} \ ,~~\L_{24} = \L_{28} = \L_{18} - m_5 \a_{36}  \nn\\ 
   _3\chi^{27}_7 &=& \{ m_{2}, m_{45}, m_{6}, -m_{12,46}, m_{12}, m_{4} \   \} \ ,~~\L_{27} = \L_{30} =\L_{22} - m_5 \a_{36} \nn\\ 
  _3\chi^{34}_7 &=& \{ m_{4}, m_{5}, m_{6}, -m_{12,46}, m_{1}, m_{2} \   \} \ ,~~\L_{34} = \L_{35} = \L_{32} - m_4 \a_{26} \nn\\  
  \nn\eea 

\eqnn{slt74}
 _4\chi^6_7 &=& \{ m_1, m_{23}, m_{5}, -m_{3,5}, m_{3}, m_{56} \   \} \ ,~~\L_6 = \L_{10} =\L_3 - m_5 \a_{45} \\ 
 _4\chi^9_7 &=& \{ m_{12}, m_{3}, m_{5}, -m_{23,5}, m_{23}, m_{56} \   \} \ ,~~\L_9 = \L_{14} = \L_5 - m_5 \a_{45} \nn\\ 
 _4\chi^{11}_7 &=& \{ m_1, m_{23}, m_{56}, -m_{3,56}, m_{3}, m_{5} \   \} \ ,~~\L_{11} = \L_{15} = \L_6 - m_6 \a_{46}  \nn\\ 
_4\chi^{12}_7 &=& \{ m_2, m_{3}, m_{5}, -m_{13,5}, m_{13}, m_{56} \   \} \ ,~~\L_{12} = \L_{16} =\L_8 - m_5 \a_{45} \nn\\ 
   _4\chi^{17}_7 &=& \{ m_{2}, m_{3}, m_{56}, -m_{13,56}, m_{13}, m_{5} \   \} \ ,~~\L_{17} = \L_{22} =\L_{13} - m_1 \a_{14} \nn\\ 
_4\chi^{18}_7 &=& \{ m_{12}, m_{3}, m_{56}, -m_{23,56}, m_{23}, m_{5} \   \} \ ,~~\L_{18} = \L_{13} =\L_{15} - m_2 \a_{24}  \nn\\ 
   _4\chi^{28}_7 &=& \{ m_{13}, m_{5}, m_{6}, -m_{23,56}, m_{2}, m_{3} \   \} \ ,~~\L_{28} = \L_{31} =\L_{24} - m_3 \a_{25}  \nn\\ 
 _4\chi^{30}_7 &=& \{ m_{23}, m_{5}, m_{6}, -m_{13,56}, m_{12}, m_{3} \   \} \ ,~~\L_{30} = \L_{33} =\L_{26} - m_5 \a_{36} \nn\\ 
 _4\chi^{32}_7 &=& \{ m_{3}, m_{5}, m_{6}, -m_{13,56}, m_{1}, m_{23} \   \} \ ,~~\L_{32} = \L_{34} =\L_{29} - m_5 \a_{36}
   \nn\eea 

 \eqnn{slt75}
  _5\chi^2_7 &=& \{ m_1, m_2, m_{34}, -m_4, m_{4}, m_6 \   \} \ ,~~\L_2 = \L_4 =\L_1 - m_4 \a_4 \\ 
 _5\chi^3_7 &=& \{ m_1, m_{23}, m_{4}, -m_{34}, m_{34}, m_6 \   \} \ ,~~\L_3 = \L_6 = \L_2 - m_3 \a_{34} \nn\\ 
    _5\chi^5_7 &=& \{ m_{12}, m_{3}, m_{4}, -m_{24}, m_{24}, m_6 \   \} \ ,~~\L_5 = \L_9 = \L_3 - m_2 \a_{24} \nn\\ 
  _5\chi^8_7 &=& \{ m_2, m_{3}, m_{4}, -m_{14}, m_{14}, m_6 \   \} \ ,~~\L_8 =\L_{12} =  \L_5 - m_1 \a_{14} \nn\\ 
  _5\chi^{15}_7 &=& \{ m_1, m_{24}, m_{6}, -m_{36}, m_{3}, m_{4} \   \} \ ,~~\L_{15} = L_{20} =\L_{10} - m_6 \a_{46}  \nn\\ 
 _5\chi^{18}_7 &=& \{ m_{12}, m_{34}, m_{6}, -m_{24,6}, m_{23}, m_{4} \   \} \ ,~~\L_{18} =\L_{24} =  \L_{15} - m_2 \a_{24}  \nn\\ 
  _5\chi^{22}_7 &=& \{ m_{2}, m_{34}, m_{6}, -m_{14,6}, m_{13}, m_{4} \   \} \ ,~~\L_{22} = \L_{27} = \L_{16} - m_6 \a_{46} \nn\\ 
_5\chi^{23}_7 &=& \{ m_{13}, m_{4}, m_{6}, -m_{24,6}, m_{2}, m_{34} \   \} \ ,~~\L_{23} = \L_{28} = \L_{18} - m_3 \a_{25}  \nn\\ 
  _5\chi^{26}_7 &=& \{ m_{23}, m_{4}, m_{6}, -m_{14,6}, m_{12}, m_{3} \   \} \ ,~~\L_{26} = \L_{30} =\L_{21} - m_6 \a_{46} \nn\\ 
   _5\chi^{29}_7 &=& \{ m_{3}, m_{4}, m_{6}, -m_{14,6}, m_{1}, m_{24} \   \} \ ,~~\L_{29} = \L_{32} = \L_{25} - m_6 \a_{46}
   \nn\eea 

   \eqnn{slt76}
  _6\chi^4_7 &=& \{ m_1, m_2, m_{34}, -m_{45}, m_{4}, m_{5} \   \} \ ,~~\L_4 = \L_7 = \L_2 - m_5 \a_{45}  \\
 _6\chi^6_7 &=& \{ m_1, m_{23}, m_{45}, -m_{35}, m_{34}, m_{5} \   \} \ ,~~\L_6 = \L_{11} = \L_3 - m_5 \a_{45} \nn\\ 
 _6\chi^9_7 &=& \{ m_{12}, m_{3}, m_{45}, -m_{25}, m_{24}, m_{5} \   \} \ ,~~\L_9 = \L_{13} = \L_5 - m_5 \a_{45} \nn\\ 
_6\chi^{10}_7 &=& \{ m_1, m_{24}, m_{5}, -m_{35}, m_{3}, m_{45} \   \} \ ,~~\L_{10} = \L_{15} = \L_6 - m_4 \a_{35}  \nn\\ 
 _6\chi^{12}_7 &=& \{ m_2, m_{3}, m_{45}, -m_{15}, m_{14}, m_{5} \   \} \ ,~~\L_{12} = \L_{17} = \L_8 - m_5 \a_{45} \nn\\ 
   _6\chi^{16}_7 &=& \{ m_2, m_{34}, m_{5}, -m_{15}, m_{13}, m_{45} \   \} \ ,~~\L_{16} = \L_{12} - m_4 \a_{35} \nn\\ 
 _6\chi^{18}_7 &=& \{ m_{12}, m_{34}, m_{5}, -m_{25}, m_{23}, m_{45} \   \} \ ,~~\L_{18} = \L_{14} =\L_{15} - m_2 \a_{24}  \nn\\ 
  _6\chi^{21}_7 &=& \{ m_{23}, m_{4}, m_{5}, -m_{15}, m_{12}, m_{35} \   \} \ ,~~\L_{21} = \L_{26} =\L_{16} - m_3 \a_{25} \nn\\ 
 _6\chi^{23}_7 &=& \{ m_{13}, m_{4}, m_{5}, -m_{25}, m_{2}, m_{35} \   \} \ ,~~\L_{23} =\L_{19} =  \L_{18} - m_3 \a_{25}  \nn\\ 
 _6\chi^{25}_7 &=& \{ m_{3}, m_{4}, m_{5}, -m_{15}, m_{1}, m_{25} \   \} \ ,~~\L_{25} = \L_{29} =\L_{21} - m_2 \a_{15}
  \nn\eea 

  \subsection{Further reductions of multiplets}

  \eqnn{slt713}
    _{13}\chi^{14}_7 &=& \{ m_{2}, m_{4}, m_{5}, -m_{2,45}, m_{2}, m_{46} \   \} \ ,~~\L_{14} = \L_{16} =\L_9 - m_4 \a_{35} \\ 
   _{13}\chi^{18}_7 &=& \{ m_{2}, m_{4}, m_{56}, -m_{2,46}, m_{2}, m_{45} \   \} \ ,~~\L_{18} = \L_{22} = \L_{15} - m_2 \a_{24}  \nn\\ 
  _{13}\chi^{23}_7 &=& \{ m_{2}, m_{4}, m_{56}, -m_{2,46}, m_{2}, m_{5} \   \} \ ,~~\L_{23} = \L_{26} = \L_{18} - m_3 \a_{25}  \nn\\ 
_{13}\chi^{24}_7 &=& \{ m_{2}, m_{45}, m_{6}, -m_{2,46}, m_{2}, m_{4} \   \} \ ,~~\L_{24} = \L_{27} = \L_{18} - m_5 \a_{36}   
      \nn\eea 

  \eqnn{slt714}
 _{14}\chi^9_7 &=& \{ m_{2}, m_{3}, m_{5}, -m_{23,5}, m_{23}, m_{56} \   \} \ ,~~\L_9 = \L_{12} =\L_5 - m_5 \a_{45} \\ 
   _{14}\chi^{18}_7 &=& \{ m_{2}, m_{3}, m_{56}, -m_{23,56}, m_{23}, m_{5} \   \} \ ,~~\L_{18} = \L_{22} = \L_{15} - m_2 \a_{24}  \nn\\ 
   _{14}\chi^{28}_7 &=& \{ m_{23}, m_{5}, m_{6}, -m_{23,56}, m_{2}, m_{3} \   \} \ ,~~\L_{28} = \L_{30} = \L_{24} - m_3 \a_{25}   
      \nn\eea 

      \eqnn{slt715}
   _{15}\chi^5_7 &=& \{ m_{2}, m_{3}, m_{4}, -m_{24}, m_{2}, m_6 \   \} \ ,~ \\ 
 _{15}\chi^9_7 &=& \{ m_{2}, m_{3}, m_{4}, -m_{23}, m_{24}, m_{6} \   \} \ ,~~\L_9 = \L_{12} =\L_5 - m_5 \a_{45} \nn\\ 
      _{15}\chi^{18}_7 &=& \{ m_{2}, m_{34}, m_{6}, -m_{24,6}, m_{23}, m_{4} \   \} \ ,~~\L_{18} = \L_{22} = \L_{15} - m_2 \a_{24}  \nn\\ 
  _{15}\chi^{23}_7 &=& \{ m_{23}, m_{4}, m_{6}, -m_{24,6}, m_{2}, m_{34} \   \} \ ,~~\L_{23} = \L_{26} = \L_{18} - m_3 \a_{25}    
        \nn\eea 

      \eqnn{slt716}
    _{16}\chi^9_7 &=& \{ m_{2}, m_{3}, m_{45}, -m_{25}, m_{24}, m_{5} \   \} \ ,~~\L_9 = \L_{12} =\L_5 - m_5 \a_{45} \\ 
      _{16}\chi^{18}_7 &=& \{ m_{2}, m_{34}, m_{5}, -m_{25}, m_{23}, m_{45} \   \} \ ,~~\L_{18} = \L_{22} = \L_{15} - m_2 \a_{24}  \nn\\ 
_{16}\chi^{19}_7 &=& \{ m_{23}, m_{4}, m_{5}, -m_{25}, m_{2}, m_{35} \   \} \ ,~~\L_{19} = \L_{21} = \L_{14} - m_3 \a_{25} \nn\\ 
    \nn\eea 

  \eqnn{slt724}
       _{24}\chi^6_7 &=& \{ m_1, m_{3}, m_{5}, -m_{3,5}, m_{3}, m_{56} \   \} \ ,~~\L_6 = \L_9 = \L_3 - m_5 \a_{45} \\ 
  _{24}\chi^{11}_7 &=& \{ m_1, m_{3}, m_{56}, -m_{3,56}, m_{3}, m_{5} \   \} \ ,~~\L_{11} = \L_{13} = \L_6 - m_6 \a_{46}  \nn\\ 
      _{24}\chi^{30}_7 &=& \{ m_{3}, m_{5}, m_{6}, -m_{1,3,56}, m_{1}, m_{3} \   \} \ ,~~\L_{30} = \L_{32} =\L_{26} - m_5 \a_{36}
 \nn  \eea

  \eqnn{slt725}
  _{25}\chi^3_7 &=& \{ m_1, m_{3}, m_{4}, -m_{34}, m_{34}, m_6 \   \} \ ,~~\L_3 = \L_3 = \L_2 - m_3 \a_{34} \\ 
         _{25}\chi^{15}_7 &=& \{ m_1, m_{34}, m_{6}, -m_{34,6}, m_{3}, m_{4} \   \} \ ,~~\L_{15} = \L_{18} = \L_{10} - m_6 \a_{46}  \nn\\ 
      _{25}\chi^{26}_7 &=& \{ m_{3}, m_{4}, m_{6}, -m_{34,6}, m_{1}, m_{3} \   \} \ ,~~\L_{26} = \L_{21} - m_6 \a_{46} \nn\\ 
  _{25}\chi^{29}_7 &=& \{ m_{3}, m_{4}, m_{6}, -m_{1,34,6}, m_{1}, m_{34} \   \} \ ,~~\L_{29} = \L_{25} - m_6 \a_{46}
 \nn   \eea

\eqnn{slt726}
       _{26}\chi^6_7 &=& \{ m_1, m_{3}, m_{45}, -m_{35}, m_{34}, m_{56} \   \} \ ,~~\L_6 = \L_9 = \L_3 - m_5 \a_{45} \\ 
  _{26}\chi^{10}_7 &=& \{ m_1, m_{34}, m_{5}, -m_{35}, m_{3}, m_{46} \   \} \ ,~~\L_{10} = \L_{14} =\L_6 - m_4 \a_{35}  \nn\\ 
 _{26}\chi^{21}_7 &=& \{ m_{3}, m_{4}, m_{5}, -m_{1,35}, m_{1}, m_{36} \   \} \ ,~~\L_{21} = \L_{25} =\L_{16} - m_3 \a_{25} \nn 
   \eea

\eqnn{slt735}
 _{35}\chi^2_7 &=& \{ m_1, m_2, m_{4}, -m_4, m_{4}, m_6 \   \} \ ,~~\L_2 = \L_1 - m_4 \a_4 \\ 
     _{35}\chi^{18}_7 &=& \{ m_{12}, m_{4}, m_{6}, -m_{2,4,6}, m_{2}, m_{4} \   \} \ ,~~\L_{18} = \L_{23} = \L_{15} - m_2 \a_{24}  \nn\\ 
 _{35}\chi^{22}_7 &=& \{ m_{2}, m_{4}, m_{6}, -m_{12,4,6}, m_{12}, m_{4} \   \} \ ,~~\L_{22} = \L_{26} = \L_{16} - m_6 \a_{46}  
   \nn\eea 

\eqnn{slt736}
    _{36}\chi^4_7 &=& \{ m_1, m_2, m_{4}, -m_{45}, m_{4}, m_{5} \   \} \ ,~~\L_4 = \L_2 - m_5 \a_{45} \\ 
   _{36}\chi^6_7 &=& \{ m_1, m_{2}, m_{45}, -m_{45}, m_{4}, m_{5} \   \} \ ,~~\L_6 = \L_3 - m_5 \a_{45} \nn\\ 
 _{36}\chi^{16}_7 &=& \{ m_2, m_{4}, m_{5}, -m_{12,45}, m_{12}, m_{45} \   \} \ ,~~\L_{16} = \L_{21} = \L_{12} - m_4 \a_{35} \nn\\ 
  _{36}\chi^{18}_7 &=& \{ m_{12}, m_{4}, m_{5}, -m_{2,45}, m_{2}, m_{45} \   \} \ ,~~\L_{18} = \L_{23} = \L_{15} - m_2 \a_{24}
    \nn\eea 

 \eqnn{slt746}
 _{46}\chi^6_7 &=& \{ m_1, m_{23}, m_{5}, -m_{3,5}, m_{3}, m_{5} \   \} \ ,~~\L_6 = \L_{10} =\L_3 - m_5 \a_{45} \\ 
  _{46}\chi^{12}_7 &=& \{ m_2, m_{3}, m_{5}, -m_{13,5}, m_{13}, m_{5} \   \} \ ,~~\L_{12} = \L_{16} =\L_8 - m_5 \a_{45} \nn\\ 
  _{46}\chi^{18}_7 &=& \{ m_{12}, m_{3}, m_{5}, -m_{23,5}, m_{23}, m_{5} \   \} \ ,~~\L_{18} = \L_{13} =\L_{15} - m_2 \a_{24}  \nn\\ 
     \nn\eea 

     \eqnn{slt713456}
       _{135}\chi^{18}_7 &=& \{ m_{2}, m_{4}, m_{6}, -m_{2,4,6}, m_{2}, m_{4} \   \} \ ,~~\L_{18} = \L_{22} = \L_{15} - m_2 \a_{24}
          \\
          _{135}\chi^{18}_7 &=& \{ m_{2}, m_{4}, m_{5}, -m_{2,45}, m_{2}, m_{45} \   \} \ ,~~\L_{18} = \L_{22} = \L_{15} - m_2 \a_{24}  \nn\\ 
  _{135}\chi^{23}_7 &=& \{ m_{2}, m_{4}, m_{5}, -m_{2,45}, m_{2}, m_{5} \   \} \ ,~~\L_{23} = \L_{26} = \L_{18} - m_3 \a_{25}  \nn\\ 
             _{135}\chi^{18}_7 &=& \{ m_{2}, m_{3}, m_{5}, -m_{23,5}, m_{23}, m_{5} \   \} \ ,~~\L_{18} = \L_{22} = \L_{15} - m_2 \a_{24}
          \nn\eea 

\section*{Conclusion and Outlook}

The ERs listed in the previous section exhaust all  representations of the algebras ~$sl(2n+1)$ for ~$n=1,2,3$.
 induced in the case ~$\cm = sl(n+1,n)$. 
 
 We noted the differences with respect to the cases $sl(2n)$ in appropriate Remarks.

Obviously, there many more cases of Langlands duals to which our approach  can be applied. This will be done in some future papers.

\newpage



$$ .\vspace{-20cm}  \includegraphics[width=15cm]  {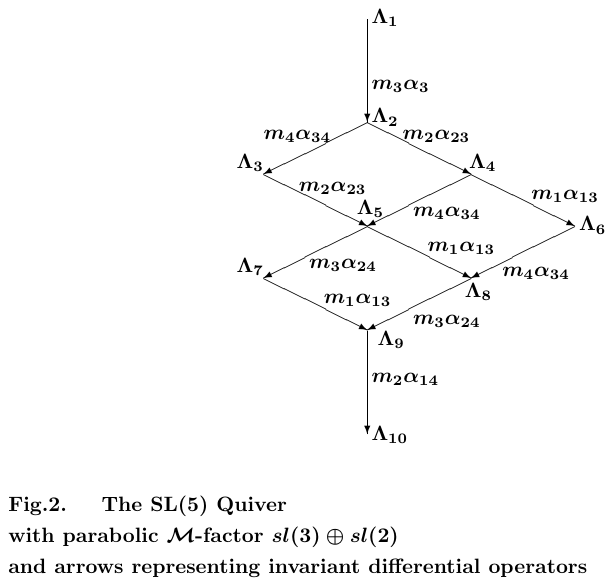}$$

\newpage



$$ .\vspace{-20cm}  \includegraphics[width=15cm]  {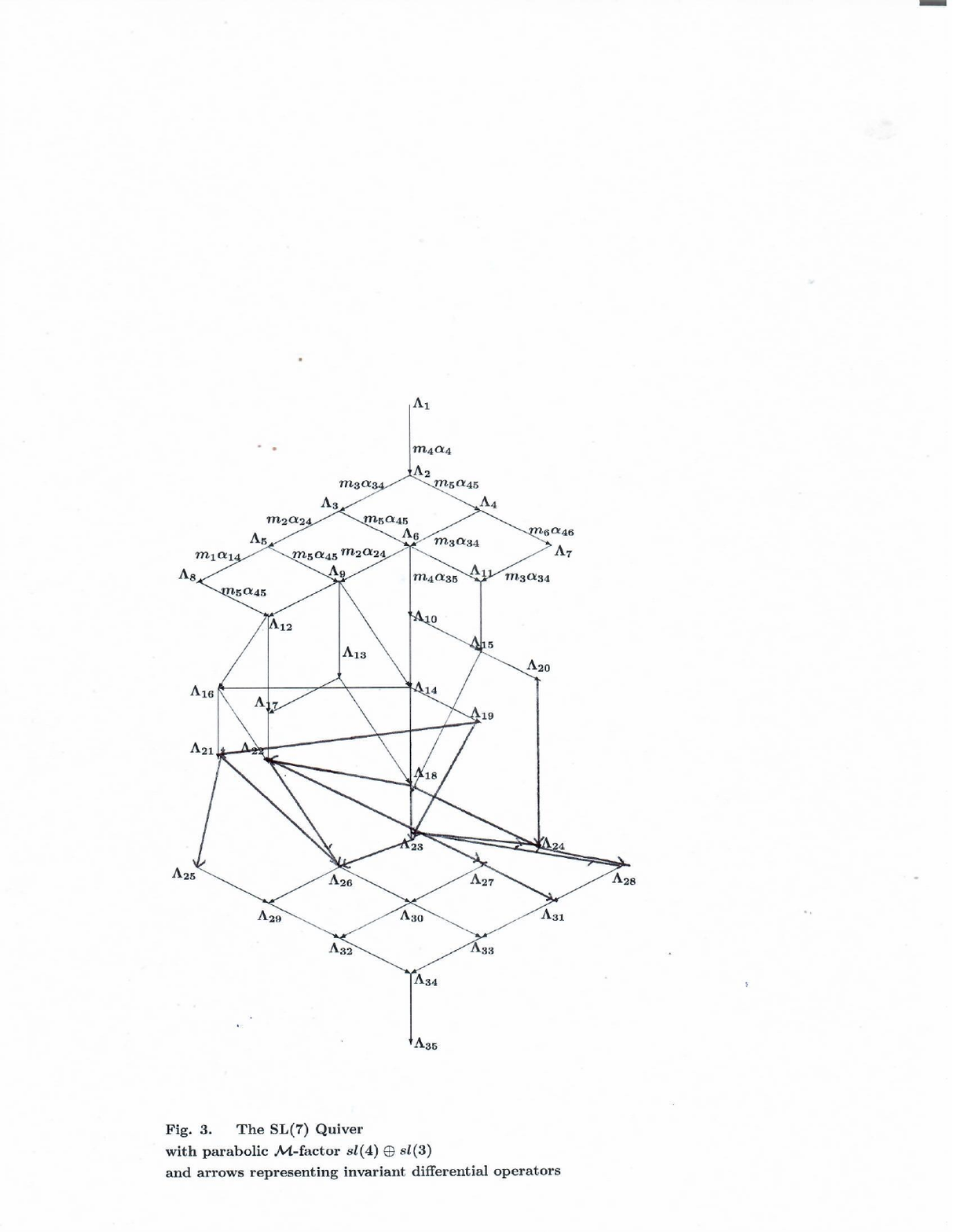}$$

\end{document}